\documentclass[12pt]{article}

\newtheorem{theorem}{Theorem}

\newtheorem{proposition}[theorem]{Proposition}

\newtheorem{corollary}[theorem]{Corollary}

\begin{document}

\author{T. Kadeishvili}

\title{Structure of $A(\infty )$-algebra and Hochschild and Harrison
cohomology}

\date{ }

\maketitle

In \cite{Kade1} for an arbitrary differential algebra $C$
(with free homology modules) in the homology algebra $H(C)$ we
have constructed a sequence of operations
$\{m_i:\otimes^iH ( C) \to H( C),\ i=3,4,...\}$,
which, together with ordinary multiplication
$m_2:H( C)\otimes H( C)\to H( C)$, turns $H( C)$
into an $A(\infty )$-algebra in the sense of Stasheff \cite{Stash}.
If a differential algebra $C$ is commutative then on $H( C)$
arises an $A(\infty )$-algebra structure of special type which
we call commutative. Particularly $A(\infty )$-algebra structure
arises on the cohomology algebra $H^*(B,\Lambda )$ of a topological
space, and a commutative $A(\infty )$-algebra structure arises  on
the rational cohomology algebra $H^*(B,Q)$. Clearly the
$A(\infty )$-algebra  $(H^*(B,\Lambda ),\{m_i\})$ carries more
information than algebra $H^*(B,\Lambda )$. Particularly cohomology
$A(\infty )$-algebra  $(H^*(B,\Lambda ),\{m_i\})$ determines
cohomology groups of the loop space $\Omega B$, and commutative 
$A(\infty )$-algebra  $(H^*(B,Q),\{m_i\})$ determines the rational
homotopy type of $B$. Naturally arises a question when these
structures are degenerate, that is when for an $A(\infty )$-algebra
$(H(C ), \{m_i\})$ the operations $m_i,\ i\geq 3$ are trivial?

In this paper we study the connection between $A(\infty )$-structures and 
Hochschild (Harrison in commutative case) cohomology of the algebra $H(C)$. 
Particularly we show that if Hochschild cohomology $Hoch^{n,2-n}(H(C),H(C))=0$ 
for $n\geq 3$ then any $A(\infty )$-algebra structure on $H(C)$ is degenerate. 
Respectively, in commutative case, any commutative $A(\infty )$-algebra 
structure is degenerate whenever Harrison cohomology 
$Harr^{n,2-n}(H(C),H(C))=0$ for $n\geq 3$

Bellow $\Lambda $ denotes a field. If $M=\sum M_q$ is a graded
$\Lambda $-module and $a\in M_p$ then $\hat{a}$ denotes $(-1)^{dim a}$.
To the permutation of elements $a\in M_p,\ b\in M_q$ corresponds
the sign $(-1)^{pq}$, this rule assignees to an arbitrary permutation
$\sigma $ of graded elements $a_1,...,a_n$ the sign denoted by $\epsilon (\sigma )$
(the Koszul sign).

\section{Products}

For an arbitrary graded module $M=\sum M_q$  the tensor coalgebra $T(M)$ is defined as
$$
T(M)=\Lambda +M+M\otimes M+...=\sum_{i=0}^{\infty }\otimes ^i M;
$$
the grading in $T(M)$ is defined by
$dim (a_1\otimes ... \otimes a_n)=\sum dima_i-n$, and the
comultipliciation $\nabla :T(M)\to T(M)\otimes T(M)$ looks as
$$
\nabla (a_1\otimes ... \otimes a_n)=\sum_{k=0}^{n }
( a_1\otimes ... \otimes a_k)\otimes (a_{k+1}\otimes ... \otimes a_n)
$$
(here the empty bracket $()$ means $1\in \Lambda $).
Iterating the comuultiplication $\nabla $ we obtain a sequence of homomorphisms
$$
\{ \nabla^i :T(M)\to \otimes^i T(M),\ i=1,2,...  \}
$$
where
$$
\nabla^1=id,  \ \nabla^2=\nabla , \ \nabla^n= (id\otimes \nabla^{n-1} )\nabla.
$$
There is also a product $\mu :T(M)\otimes T(M)\to T(M)$ which, together
with $\nabla $, determines on $T(M)$ a structure of Hopf algebra, this is
the shuffle product defined by Eilenberg and MacLane in\cite{Eil-Mac}.
This product is defined as
$$
\mu ((a_1\otimes ... \otimes a_n)\otimes (a_{n+1}\otimes ... \otimes a_{n+m}))=
\sum \epsilon (\sigma ) a_{\sigma (1)}\otimes ... \otimes a_{\sigma (n+m)},
$$
where summation is taken over all permutations of the set $(1,2,...,n+m)$
which satisfy the condition: $i<j$ if $1\leq \sigma (i)<\sigma (j)\leq n$
or  $n+1\leq \sigma (i)<\sigma (j)\leq n+m$.

In \cite{comm} it is shown that this product is uniquely  characterized
by the following axioms:

(1) $\mu $ turns $T(M)$ into a Hopf algebra;

(2) $p_0(\mu (a_+\otimes b_+))= p_1(\mu (a_+\otimes b_+))=0$
for arbitrary $a_+,b_+\in T(M)_+$ (here $T(M)_+=\sum_{n\geq 1} \otimes ^nM$
and $p_i:T(M)\to \otimes ^iM$ is the clear projection).

Let us denote
$$
\begin{array}{c}
Sh^n(M)=\sum_{k=1}^{n-1}\mu ((\otimes^kM)\otimes (\otimes^{n-k}M))\subset \otimes^nM;\\
Sh(M)=\sum Sh^n(M);\ Ch^n(M)=\otimes^nM/Sh^n(M);\ Ch(M)=\sum Ch^n(M).
\end{array}
$$
It is clear that
$$
\begin{array}{c}
Sh(M)=\mu(T(M)_+\otimes T(M)_+);\ Ch(M)=T(M)/Sh(M);\\
Sh^0(M)=Sh^1(M)=0;\ Ch^0(M)=\Lambda; Ch^1(M)=M.
\end{array}
$$

Suppose now that together with $T(M)$ there is given a graded
$\Lambda $-module $N=\sum N_q$. Denote by $Hom^k(\otimes^nM,N)$ a set
of $\Lambda $-homomorphisms $f:\otimes^nM\to N$ of degree $k$
(that is $f(a_1\otimes ... \otimes a_n)\subset M_{q+k}$ where
$q=\sum dim a_i$). We also use the notation
$$
Hom(\otimes^nM,N)=\sum_{k}Hom^k(\otimes^nM,N),\ Hom(T(M),N)=\sum_nHom(\otimes^nM,N).
$$
In $Hom(T(M),M)$ Gerstenhaber \cite{Gerst} has introduced a product
which we denote by $\smile_1$, defined as follows. For
$f\in Hom^k(\otimes^mM,N),\ g\in Hom^t(\otimes^nM,N)$ the product
$f\smile_1g\in Hom^{k+t}(\otimes^{m+n-1}M,N)$ looks as
$$
f\smile_1g(a_1\otimes ... \otimes a_{m+n-1})=
\sum_{k=0}^{m-1}\pm f(a_1\otimes ...\otimes a_k\otimes g(a_{k+1}\otimes ...
\otimes a_{k+n})\otimes ... \otimes a_{m+n-1}).
$$

We shall use also more general product of an element $f\in Hom(T(M),M)$
and a sequence $(g_1,...,g_k),\ g_i\in Hom(T(M),M)$ defined as
$$
\begin{array}{l}
f\smile_1(g_1,...,g_k)(a_1\otimes ... \otimes a_n)=\\
\sum \pm f(a_1\otimes ...\otimes a_{k_1}\otimes g(a_{k_1+1}\otimes ...
\otimes a_{k_1+p_1})\otimes a_{k_1+p_1+1}\otimes ... \otimes a_{k_1+t_1}\otimes\\
g_2(a_{k_1+t_1+1}\otimes ...\otimes a_{k_1+t_1+p_2})\otimes ...\\
\otimes a_s\otimes
g_k(a_{s+1}\otimes ...\otimes a_{s+p_k})\otimes a_{s+p_k+1}\otimes ...
\otimes a_n).
\end{array}
$$

The product $f\smile_1 g$ can be written as
$f\smile_1g=f(id\otimes g\otimes id)\nabla^3$. Similarly
$ f\smile_1(g_1,...,g_k)=f(id\otimes g_1\otimes id\otimes g_2\otimes ...
\otimes id\otimes g_k\otimes id)\nabla^{2k+1}.$
We remark also that these products are defined if $f\in Hom(T(M),N)$ and
$g_i\in Hom(T(M),M)$.

The product $f\smile_1g$ is not  associative generally but easy to see that
$$
f\smile_1(g\smile_1h)- (f\smile_1g)\smile_1h=
f\smile_1(g,h)- f\smile_1(h,g).
$$

Let us consider $Hom(Ch(M),M)=\sum_{n} Hom(Ch^n(M),M)$.
Clearly
$$
Hom(Ch^n(M),M)=\{f\in Hom(\otimes^n(M),M),\ f|Sh^n(M)=0\},
$$
hence $Hom(Ch^n(M),M)$ is also graded by degree of homomorphisms, i.e.
$$
 Hom(Ch^n(M),M)=\sum_k Hom^k(Ch^n(M),M)
$$
where
$$
 Hom^k(Ch^n(M),M)= \{f\in Hom^k(\otimes^nM,M),\ f|Sh^n(M)=0\},
$$
thus
$ Hom(Ch(M),M)\subset Hom(T(M),M)$. It is possible to show that if
$f,g\in Hom(Ch(M),M)$ then $f\smile_1g\in Hom(Ch(M),M)$. Moreover,
$f\smile_1(g,...,g)\in Hom(Ch(M),M)$ too.

\section{Hochshild cohomology}

Let $A$ be a graded algebra and $M$ be a graded bimodule over $A$.
Hochschild cochain complex  is defined as
$$
C^*(A,M)=\sum C^n(A,M),\  \ C^n(A,M)=Hom(\otimes^nA,M),
$$
the coboundary operator $\delta: C^n(A,M)\to C^{n+1}(A,M)$ is given by
$$
\begin{array}{c}
\delta f(a_1\otimes ... \otimes a_{n+1})=
a_1f((a_2\otimes ...
\otimes a_{n+1})+\\
\sum_{k}\pm f(a_1\otimes ... \otimes a_ka_{k+1}\otimes  ...
\otimes a_{n+1})\pm f(a_1\otimes ... \otimes a_{n})a_{n+1}.
\end{array}
$$
Hochschild cohomology of $A$ with coefficients in $M$ is defined as homology of this cochain complex and is denoted by $Hoch^*(A,M)$. Since $A$
and $M$ are graded, each $C^n(A,M)$ is graded too:
$C^n(A,M)=\sum_{k}C^{n,k}(A,M)$ where $C^{n,k}(A,M)=Hom^k(\otimes^nA,M)$.
It is easy to see that $\delta: C^{n,k}(A,M)\to C^{n+1,k}(A,M)$ , so
$(C^{*,k}(A,M),\delta) $ is a direct summand in
$(C^*(A,M),\delta)$, thus Hochschild cohomology in this case is bigraded: 
$Hoch^n(A,M)=\sum_{k} Hoch^{n,k}(A,M)$
where $ Hoch^{n,k}(A,M)$ is the $n$-th homology module of
$ (C^{*,k}(A,M),\delta )$. 

Instead of $M$ we can take the algebra $A$ itself. The complex
$ C^{*,*}(A,A)$ is a differential algebra with respect to the
following product:
for $f\in C^{m,k}(A,A)$ and $g\in C^{n,t}(A,A)$ the product
$f\smile g\in C^{m+n,k+t}(A,A)$ is defined by
$$
f\smile g(a_1\otimes ... \otimes a_{m+n})=f(a_1\otimes ...
\otimes a_{m})\cdot g (a_{n+1}\otimes ... \otimes a_{m+n}).
$$
Besides this product in $ C^{*,*}(A,A)$ we have also the product
$f\smile_1g$ (see the previous section). In \cite{Gerst} it is shown that these
products satisfy the standard conditions
\begin{equation}
\label{steen0}
\delta (f\smile g)=\delta f\smile g \pm f\smile \delta g;
\end{equation}
\begin{equation}
\label{steen1}
\delta (f\smile_1 g)=\delta f\smile_1 g \pm f\smile_1 \delta g
\pm f\smile g\pm g\smile f.
\end{equation}

\section{Harrison cohomology}

Suppose now $A$ is a commutative graded algebra and $M$ 
is a module over $A$. The Harrison cochain complex $\bar{C}^*(A,M)$
is defined as a subcomplex of the Hochschild complex
$$
\bar{C}^*(A,M)=\{f\in C^*(A,M)|\ f|Sh (A)=0\},
$$ 
i.e.
$\bar{C}^*(A,M)=Hom(Ch(A),M)=\sum_{n}Hom(Ch^n(A),M).$
In \cite{Harr} it is shown that $\bar{C}^*(A,M)$ is closed
with respect to the differential $\delta $ (that is if $ f|Sh (A)=0$
then $\delta f|Sh (A)=0$). Harrison cohomology of a commutative algebra
$A$ with coefficients in an $A$-module $M$ is defined as homology of the
cochain complex $(\bar{C}^*(A,M),\delta )$. These cohomologies we denote by
$Harr^*(A,M)$. As above each module $\bar{C}^n(A,M)=Hom(Ch^n(A),M)$ is
graded by degrees of homomorphisms:
$\bar{C}^n(A,M)=\sum_{k}\bar{C}^{n,k}(A,M)$ where
$\bar{C}^{n,k}(A,M)=Hom^k(Ch^n(A),M)$. Besides
$\delta: \bar{C}^{n,k}(A,M)\to \bar{C}^{n+1,k}(A,M)$ , so
$(\bar{C}^{*,k}(A,M),\delta) $ is a direct summand in the Harrison complex
$(\bar{C}^n(A,M),\delta)$. Thus
$Harr^n(A,M)=
\sum_{k} Harr^{n,k}(A,M)$ where $ Harr^{n,k}(A,M)$
is the $n$-th homology module of
$ (\bar{C}^{*,k}(A,M),\delta )$. So Harrison cohomology in is
bigraded for graded $A$ and $M$.

As it is mentioned in  the section 1 the Harrison complex $(\bar{C}^{*,k}(A,M)$ is
closed with respect to products $f\smile_1g$ and $f\smile_1(g,...,g)$.
The formulae (\ref{steen0}) and (\ref{steen1}) are valid in this subcomplex too.

\section{Twisting cochains in the Hochschild and Harrison complexes}

In \cite{Berik1}, \cite{Berik2} N. Berikashvilihas has defined a functor $D$ from
the category of differential algebras to the category of pointed sets, which
have applications in the homology theory of fibrations. We recall shortly
its definition. Let $(C,d)$ be a differential graded algebra. A twisting
cochain is defined as an element $a=a^2+a^3+...,\ a^i\in C^i$ satisfying
the condition $da=\pm a\cdot a$. Let $Tw( C)$ be the set of all twisting cochains.
In this set by Berikashvili was introduced the following equivalence relation:
$a\sim  a´$ if there exists an element $p=p^1+p^2+...,\ p^i\in C^i$ such that
$$
a-a`=p\cdot a \pm a´\cdot p \pm dp.
$$
The factorset of the set $Tw( C)$ by this equivalence is denoted by $D( C)$.
We are going to introduce the similar definition in the Hochschild and
Harrison complexes but with respect to $\smile_1$ product. Note that Hochschild and Hasrrison complexes are
not differential algebras with respect to the product $f\smile_1 g$,
besides this product is not associative, hence in order to define the
functor $D$ some modification is needed.

Let us define a twisting cochain in the Hochschild complex $C^{*,*}(A,A)$
as an element $a=a^{3,-1}+ a^{4,-2}+...+ a^{i,2-i}+...,\ a^{i,2-i}\in C^{i,2-i}(A,A)$
satisfying the condition $\delta a=a\smile_1a$. Let $Tw( A,A)$ be the set of all
twisting cochains. Now we introduce in this set the following equivalence
relation: $a\sim  a'$ if there exists an element
$p=p^{2,-1}+p^{3,-2}+...+p^{i,1-i}+...,\ ,\ p^{i,1-i}\in C^{i,1-i}(A,A)$
such that
\begin{equation}
\label{eqv}
a-a'=\delta p\pm p\smile_1a\pm a'\smile_1p\pm a'\smile_1(p,p)\pm a'\smile_1(p,p,p)\pm ...
\end{equation}
(the sum is finite in each dimension). This is an equivalence relation;
we denote by $D(A,A)$ the factorset $Tw(A,A)/\sim  $. The set $D(A,A)$ is a
pointed set: a distinguished point is the class of $a=0$, which we denote by $0\in D(A,A)$.

There is a possibility to perturb twisting cochains without changing
their equivalence classes in $D(A,A)$. Indeed, let $A\in Tw(A,A)$ and
$ p\in C^{n,1-n}(A,A)$ be an arbitrary cochain, then there exists a twisting
cochain $\bar{a}\in Tw(A,A)$ such that $a^i=\bar{a}^i$ for $i\leq n$,
$ \bar{a}^{n+1}=a^{n+1}+\delta p$ and $\bar{a}\sim  a$.
The twisting cochain $\bar{a}$ can be solved inductively from the
equation (\ref{eqv}).

\begin{theorem}\label{hoch3} If $Hoch^{n,2-n}(A,A)=0$ for $n\geq 3$ then $D(A,A)=0$.
\end{theorem}

\noindent {\bf Proof.} We have to show that in this case an arbitrary twisting cochain
is equivalent to zero. From the equality $\delta a=a\smile_1 a$
in dimension $n=3$ we obtain $\delta a^3=0$ that is $a^3\in C^{3,-1}(A,A)$
is a cocycle. Since $Hoch^{3,-1}(A,A)=0$ there exists $p^{2,-1}\in C^{2,-1}(A,A)$
such that $a^3=\delta p^{2,-1}$. Perturbing our twisting cochain $a=a^3+a^4+... $
by $p^{2,-1}$ we we obtain new twisting cochain $\bar{a}=\bar{a}^3+\bar{a}^4+...$
equivalent to $a$ and with $\bar{a}^3=0$. Now the component $\bar{a}^4$
becomes a cocycle, which, which can be killed using $Hoch^{4,-2}(A,A)=0$ etc.
This completes the proof.

Now turn to the Harrison complex $\bar{C}(A,A)$ of a commutative algebra
$A$. Since $\bar{C}(A,A)$ is closed with respect to $\smile_1$ product,
here we also can define twisting cochains and the set $\bar{D}(A,A)$
and prove the similar

\begin{theorem} If $Harr^{n,2-n}(A,A)=0$ for $n\geq 3$ then $D(A,A)=0$.
\end{theorem}

\section{ Structure of $A(\infty )$-algebra and Hochschild and Harrison cohomology }

$A(\infty )$-algebra was defined by Stasheff in\cite{Stash}.
It is a graded $\Lambda $-module $M=\sum M_q$ equipped with a
sequence of operations - $\Lambda $-homomorphisms
$\{m_i:\otimes^iM\to M,\ i=1,2,...\}$ satisfying the following conditions
\begin{equation}
\label{dim}
m_i(\otimes^iM)_q\subset M_{q-i+2},\ i.e.\ \ degm_i=2-i;
\end{equation}
\begin{equation}
\label{sta}
\sum_{j=1}^n\sum_{k=0}^{n-j} \pm m_i(a_1\otimes ...\otimes a_k\otimes 
m_j(a_{k+1}\otimes ... \otimes a_{k+j})\otimes ...\otimes a_n)=0.
\end{equation}

A morphism of $A(\infty )$-algebras $f: (M,\{m_i\})\to (M',\{m'_i\})$ is a sequence of homomorphisms $\{f_i :\otimes^iM\to M',\ i=1,2,...\}$
satisfying the following conditions
\begin{equation}
\label{dimf}
f_i(\otimes^iM)_q\subset M'_{q-i+1},\ \ i.e.\ \ degf_i=1-i;
\end{equation}
\begin{equation}
\label{mor}
\begin{array}{c}
\sum_{j=1}^n\sum_{k=0}^{n-j} \pm f_i(a_1\otimes ...\otimes a_k\otimes 
m_j(a_{k+1}\otimes ... \otimes a_{k+j})\otimes ...\otimes a_n)=\\
\sum_{t=1}^n \sum_{k_1+...+k_t=n} \pm m'_t(f_{k_1}(a_1\otimes ...\otimes a_{k_1})\otimes ...\otimes f_{k_t}(a_{n-k_t+1}\otimes ...\otimes a_{n})).
\end{array}
\end{equation}
The obtained category is denoted by $A(\infty )$.  

For an arbitrary $A(\infty )$-algebra $(M,\{m_i\})$ the sequence of 
operations $\{m_i\}$ defines on the tensor coalgebra $T(M)$ a 
differential $d:T(M)\to T(M)$ given by
$$
d(a_1\otimes ...\otimes a_n)=\sum_{k,j}\pm a_1\otimes ...\otimes a_k\otimes 
m_j(a_{k+1}\otimes ... \otimes a_{k+j})\otimes ...\otimes a_n
$$
which fits with the coproduct $\nabla :T(M)\to T(M)\otimes T(M)$, i.e. turns $T(M)$ into a differential coalgebra. This differential coalgebra $(T(M),d)$ is called $\tilde{B}$-construction of $A(\infty )$-algebra $(M,\{m_i\})$ and is denoted by $\tilde{B}(M,\{m_i\})$ (\cite{Stash}).

An arbitrary morphism of $A(\infty )$-algebras $\{f_i\}:(M,\{m_i\})\to (M',\{m'_i\})$ induces a $DG$-coalgebra morphism $\tilde{B}(f):\tilde{B}(M,\{m_i\})\to \tilde{B}(M',\{m'_i\})$  by
$$
\tilde{B}(f)(a_1\otimes ...\otimes a_n)=
\sum_{t=1}^n \sum_{k_1+...+k_t=n}f_{k_1}(a_1\otimes ...\otimes a_{k_1})\otimes ...\otimes f_{k_t}(a_{n-k_t+1}\otimes ...\otimes a_{n}).
$$
Thus $\tilde{B}$ is a functor from the category of $A(\infty )$-algebras to the category of $DG$-coalgebras.

We are interested in $A(\infty )$-algebras of type $(M,\{m_1=0,m_2,m_3,...\})$, i.e. with $m_1=0$. The full subcategory of $A(\infty )$ which objects are such $A(\infty )$-algebras we denote by $A^0(\infty )$.

Now let $(M,\mu )$ be a graded associative algebra with multiplication $\mu :M\otimes M\to M$. Consider all possible $A(\infty )$ structures $\{m_i\}$ on $M$ with $m_1=0,\ m_2=\mu $. Two such structures we call equivalent  if there exists a morphism of $A(\infty )$-algebras $\{p_i\}:(M, \{m_i\} )\to (M,\{m'_i\})$ for which the first component $p_1$ is the identity map (one can show that such morphisms are isomorphisms in the category $A^0(\infty )$). The obtained factorset we denote by $(M,\mu )(\infty )$. A trivial $A^0(\infty )$ structure we define as $\{m_i\}$ with $m_{i>2}=0$. It's class we denote as $0\in (M,\mu )(\infty )$.

\begin{proposition} The sets $(M,\mu )(\infty )$ and $D(M,M)$ are bijective.
\end{proposition}

\noindent {\bf Proof.} Let $(M, \{m_i\})$ be an $A(\infty )$-algebra 
with with $m_1=0,\ m_2=\mu $; we denote
$m=m_3+m_4+...\ $. Each operation $m_i:\otimes^iM\to M$ can be interpreted as 
a Hochschild cochain from $C^i(M,M)$. It is easy to mention that the condition 
(\ref{sta}) means exactly $\delta m=m\smile_1 m$, i.e. $m\in Tw(M,M)$. 
Conversely, each twisting cochain $m=m_3+m_4+...\in Tw(M,M)$ defines on 
$m$ an $A^0(\infty )$-algebra structure $(M,\{m_1=0,m_2=\mu, m_3,m_4,...\})$. 
It remains to show that two $A^0(\infty )$-structures $\{m_i\} $ 
and $\{m'_i\} $ are equivalent if and only if the twisting cochains 
$m$ and $m'$ are equivalent. Indeed, if $A^0(\infty )$-structures 
$\{m_i\} $ and $\{m'_i\} $ are equivalent, i.e. there exists  an 
isomorphism of 
$A^0(\infty )$-algebras $\{id,p_2,p_3,...\}: (M, \{m_i\})\to (M, \{m'_i\})$, 
then the cochain $p=p_2+p_3+...$ realizes the equivalence of twisting cochains 
$m$ and $m'$, since the condition
(\ref{mor}) is equivalent to (\ref{eqv}). Conversely, if $p=p_2+p_3+...$ 
realizes the equivalence of twisting cochains $m$ and $m'$, then 
$\{id,p_2,p_3,...\}: (M, \{m_i\})\to (M, \{m'_i\})$ is the needed isomorphism.

From this proposition and the theorem \ref{hoch3} follows the

\begin{corollary} If for a graded algebra $(M,\mu )$ all $Hoch^{n,2-n}(M,M)=0$ 
for $n\geq3$ then any $A^0(\infty )$-algebra structure $\{m_i\}$ on $M$ 
(with $m_1=0,\ m_2=\mu $) is equivalent to trivial one.
\end{corollary}

Now we turn to the commutative case.

An $A(\infty )$-algebra we call commutative if the sequence of operations 
$\{m_i\}$ apart of the conditions (\ref{dim}) and (\ref{sta}) satisfies 
$m_i|Sh^i(M)=0$. In this case the differential $d:T(M)\to T(M)$ defined 
by $\{m_i\}$ fits with shuffle product, so $\tilde{B}(M,\{m_i\})$ 
becomes a DG-Hopf algebra. A morphism of commutative
$A(\infty )$-algebras $\{f_i\}:(M,\{m_i\})\to (M',\{m'_i\})$
we define as a morphism of
$A(\infty )$-algebras which, apart of the conditions (\ref{dimf}) and 
(\ref{mor}) satisfies $f_i|Sh^i(M)=0$. In this case
$\tilde{B}(f):\tilde{B}(M,\{m_i\})\to \tilde{B}(M',\{m'_i\})$
becomes a map of DG-Hopf algebras (see \cite{comm}). The condition 
$m_i|Sh^i(M)=0$ for $i=2$ means that $m_2:M\otimes M\to M$ is commutative 
and all the operations $m_i:\otimes^iM\to M$ are from the Harrison 
subcomplex $\bar{C}(M,M)$.

Now let $(M,\mu )$ be a graded commutative algebra with multiplication 
$\mu :M\otimes M\to M$. Consider all possible commutative $A(\infty )$ 
structures $\{m_i\}$ on $M$ with $m_1=0,\ m_2=\mu $; two such structures 
we call equivalent  if there exists a morphism of commutative $A(\infty )$-algebras $\{p_1=id,p_2,p_3,...\}:(M, \{m_i\} )\to (M,\{m'_i\})$ 
(which, as above, is an isomorphism). The obtained factorset we denote by 
$(M,\mu )(\infty )_c$. Exactly as above we obtain the

\begin{proposition} The sets $(M,\mu )(\infty )_c$ and $\bar{D}(M,M)$ 
are bijective.
\end{proposition}

\begin{corollary} If for a graded commutative algebra $(M,\mu )$ all 
$Harr^{n,2-n}(M,M)=0$ for $n\geq3$ then any commutative $A^0(\infty )$-algebra structure $\{m_i\}$ on $M$ (with $m_1=0,\ m_2=\mu $)is equivalent to trivial one.
\end{corollary}


\begin{thebibliography}{Harr}

\bibitem{Harr} M. Barr, Harrison homology, Hochschild homology and tripples, 
Journal of algebra, 8 (1968), 314-323.



\bibitem{Berik1}
N. Berikashvili, On differentials of spectral sequences, Bull. Georg. Acad. Sci., 
51, 1 (1968), 9-14.

\bibitem{Berik2}
N. Berikashvili, On the homology theory of spaces, Bull. Georg. Acad. Sci., 86, 
3 (1977), 529-532.

\bibitem{Eil-Mac} S. Eilenberg, S. MacLane, On the groups $H(\pi , n)$, 1, 
Ann. of Math., 58 (1953), 55-106.

\bibitem{Gerst} M. Gerstenhaber, The homology structure of an associative ring, 
Ann. of Math., 4 (1963), 267-288.

\bibitem{Kade1} T. Kadeishvili, On the homology theory of fibre spaces, 
Russian Mathematics Surveys, 6 (1980), 231-238.

\bibitem{comm} T. Kadeishvili, On the category of differential coalgebras and the 
category of $A(\infty )$-algebras, Proc. of A. Razmadze Math. Inst., 77 (1985), 50-70.



\bibitem{Stash} J. D. Stasheff, Homotopy associativity of H-spaces, TAMS, 108, 
2 (1963), 275-313.
 



\end{thebibliography}
 \end{document}